\documentclass[12pt,A4]{amsart}
\usepackage[latin1]{inputenc}
\usepackage{amssymb}
\usepackage{amsmath}
\usepackage{graphicx}

\def\R{{\Bbb R}}

\def\g{{\mathfrak g}}

\def\fracl#1#2{\frac{\displaystyle #1}{\displaystyle #2}}

\def\wt#1{\widetilde #1}

 \theoremstyle{definition}

\theoremstyle{plain}

\newtheorem{prop}{Proposition}

\newtheorem{lm}{Lemma}

\theoremstyle{remark}

\newtheorem{rem}{Remark}

\begin{document}

\title{A remark on Ricci flow of left invariant metrics}
\author{J.R. Arteaga B.}%
\address{Departamento de Matem\'aticas, Universidad de los Andes, Bogot\'a D.C., Colombia}%
\email{jarteaga@uniandes.edu.co}%
\author{M.A. Malakhaltsev}%
\address{Department of Mathematics, Kazan State University, Kremlevskaya, 18, Kazan: 420008, Russia}%
\email{Mikhail.Malakhaltsev@ksu.ru}%

\thanks{}%

\dedicatory{}%
\commby{}%
\begin{abstract}
We prove that  the Ricci flow equation for left invariant metrics on Lie
groups reduces to a first order ordinary differential equation for a
map $Q : (-a,a) \to UT$, where $UT$ is  the group of upper triangular
matrices. We decompose the matrix $R_{ij}$ of Ricci tensor
coordinates with respect to an orthonormal frame field $E_{i}$ into a sum 
$\overset{1}{R}_{ij} + \overset{2}{R}_{ij} + 
\overset{3}{R}_{ij} + \overset{4}{R}_{ij}$ such that,  
for any $E_{i'} = U^i_{i'} E_i$ with $||U^i_{i'}|| \in
O(n)$, $\overset{\alpha}{R}_{i'j'} = 
U_{i'}^i \overset{\alpha}{R}_{ij} U^j_{j'}$.
This allows us to specify several cases when the
differential equation can be simplified. 
As an example we consider three-dimensional unimodular Lie groups.
\end{abstract}
\date{July 17 2005}%
\maketitle

\section*{Introduction}
\noindent 
Let $g(t)$ be a smooth one-parameter family of metrics on  a manifold
$M$.
The Ricci flow equation defined by Richard Hamilton  in 1982 \cite{Hamilton} is:
\begin{equation}\label{rf}
\fracl{\partial}{\partial t} g_{ij}(t)=-2R_{ij}
\end{equation}
This equation can be considered as generalization of the heat equation:
\[
\fracl{\partial T}{\partial t}=\lambda \Delta T=\lambda( T_{xx} +
T_{yy} + T_{zz})
\]
where $T$ is the temperature and $\lambda$ is the thermal
diffusion.

Note that  the Ricci flow defined by the equation
(\ref{rf}) is known as an unnormalized Ricci flow. 
When the manifold $M$ is closed one can define a normalized Ricci flow 
 given by the equation
\[
\fracl{\partial}{\partial t}\widetilde{g}_{ij}=
-2R(\widetilde{g})_{ij}+\fracl{2}{n}r\widetilde{g}_{ij}
\]
where $r=\fracl{1}{Vol(M)}\int_M R(\widetilde{g})dV_{\widetilde{g}}$
\cite{Sesun}.

The Ricci flow equation was defined to study Thurston's
geometrization conjecture (for relations between the
geometrization conjecture and the Ricci flow see \cite{Chow}). 
In this framework one of the main problems is the
stability of solution, that is, 
to find  if a solution $g(t)$ of  (\ref{rf}) is stable in 
time, that is, to study whether  $g(t)$ converges as $t \to \infty$.


Now let us give several elementary examples. 

\subsection*{Example 1.}\label{example1}
Consider the sphere of radius $r$ in the Euclidean $(n+1)$-space.
The metric tensor is
\[
g_{ij}=\rho\widetilde{g}_{ij}
\]
where $\widetilde{g}$ is the standard metric of the unit sphere $S^n$,
and $\rho = r^2$. Then
\[
R_{ij}=(n-1)\widetilde{g}_{ij},
\]
hence $R_{ij}$  is independent of $r$. The Ricci flow equation (\ref{rf})
reduces to
\[
\fracl{\partial \rho}{\partial t}=-2(n-1)
\]
with solution
\[
\rho(t)=\rho(0)-2(n-1)t
\]
This example  shows that the flow is a family of spheres which
change the radius starting at the radius $r(0)$. This solution
collapses at the finite time $t=\fracl{\rho(0)}{2(n-1)}$

Using certain ideas of Hamilton such as
rescaling the metric with the volume remaining  constant, G.\,Perelman
showed that such collapse can be eliminated by performing a kind of surgery on
the manifold   and the solution converges or admits  Thurston's
decomposition
(for detailed information we refer the reader to  \cite{Hamilton}, 
\cite{Perelman}, and  \cite{Milnor}). 

\subsection*{Example 2.}\label{example2}
In Example 1, the sphere $S^n$ is an Einstein manifold, i.e.
the metric $\wt g$ satisfies the equation $Rc(\widetilde{g})=\lambda
\widetilde{g}$ 
for some $\lambda \in \mathbb{R}$. In general, if the initial metric
$\widetilde{g}$ is Einstein, then we can define the 1-parametric
family
\begin{equation}\label{ejemplo2}
g(t)=e^{-2\alpha t}\widetilde{g}, \quad \alpha=const.
\end{equation}
which is a solution of Ricci flow equation(\ref{rf}) with
$g(0)=\widetilde{g}$.

\subsection*{Example 3.}\label{example3}
Here we give  an example of solution $g(t)$ of (\ref{rf}) which is
not a family of Einstein metrics.
Let $M = \R^2$ and set 
\[
g(x,y)=\left(%
\begin{array}{cc}
  f(x) & 0 \\
  0 & \mu(x) \\
\end{array}%
\right)
\]
Then we have 
\begin{equation*}
  \begin{array}{l}
\Gamma^1_{11}=\fracl{1}{2}g^{11} f_x, \Gamma^1_{12}=0,
\Gamma^1_{22}=\fracl{1}{2}g^{11}\mu_x, \Gamma^2_{11}=0,
\Gamma^2_{12}=\fracl{1}{2}g^{22}\mu_x,
\Gamma^2_{22}=0,\\
R_{12,12}=\fracl{1}{4}\fracl{f_x}{f}\mu_x-\fracl{1}{2}\mu_{xx}+\fracl{1}{4}\fracl{\mu_x^2}{\mu}\\
R_{11}=-\fracl{1}{\mu}\left(\fracl{1}{4}\fracl{f_x}{f}\mu_x-\fracl{1}{2}\mu_{xx}+\fracl{1}{4}\fracl{\mu_x^2}{\mu}\right)\\
R_{22}=-\fracl{1}{f}\left(\fracl{1}{4}\fracl{f_x}{f}\mu_x-\fracl{1}{2}\mu_{xx}+\fracl{1}{4}\fracl{\mu_x^2}{\mu}\right)\\
R_{12}=R_{21}=0.
  \end{array}
\end{equation*}

The Ricci flow equation (\ref{rf}) reduces to 
$f_t=-2R_{11}$ and
$\mu_t=-2R_{22}$, hence follows that 
\begin{itemize}
\item $f_t=\fracl{f}{\mu}\mu_t$, which imply that $f=\lambda\mu$, and
\item $\mu \mu_{xx}-\mu_x^2-2\mu^2\mu_t=0$
\end{itemize}

In the second equation we set  $\mu(x,t)=X(x)T(t)$,
then 
\begin{equation}\label{sol3}
\mu(x,t)=\left(\fracl{A}{2}t+B\right)\left(\fracl{-1+\tanh(\fracl{x+D}{2C})^2}{2AC^2}\right)
\end{equation}
The family of metrics is:
\begin{equation}\label{metr}
g(t,x)=\left[%
\begin{array}{cc}
  f(t,x) & 0 \\
  0 & \mu(t,x) \\
\end{array}%
\right] ,\quad \textmd{where} \quad f(t,x)=\lambda \mu(t,x)
\end{equation}
Substituting this to  (\ref{rf}),  we finally get  the
solution:
\begin{equation}\label{example3}
g(t,x)= \left[
\begin {array}{cc}
2(t - A )( -1 + \tanh( x - B )^2 ) & 0 \\%
0& (t - A )( -1 + \tanh( x- B )^2 )\\
\end {array} \right]
\end{equation}
which has the Ricci tensor:
\begin{equation}\label{ricci-example3}
Rc(g(t,x))= -\left[
\begin {array}{cc}
-1 + \tanh( x - B )^2  & 0 \\%
0& \fracl{ -1 + \tanh( x- B )^2 }{2}\\
\end {array} \right]
\end{equation}
Thus the metrics $g(t)$ are not Einstein.

In the present paper we consider the  Ricci flow equation (\ref{rf}) for left
invariant metrics on Lie groups. We prove that  the Ricci flow
equation for these  metrics  reduces to a first order ordinary differential equation for a
map $Q : (-a,a) \to UT$, where $UT$ is  the group of upper triangular
matrices. We decompose the matrix $R_{ij}$ of Ricci tensor
coordinates with respect to an orthonormal frame field $E_{i}$ into a sum 
$\overset{1}{R}_{ij} + \overset{2}{R}_{ij} + 
\overset{3}{R}_{ij} + \overset{4}{R}_{ij}$ such that,  
for any $E_{i'} = U^i_{i'} E_i$ with $||U^i_{i'}|| \in
O(n)$, $\overset{\alpha}{R}_{i'j'} = 
U_{i'}^i \overset{\alpha}{R}_{ij} U^j_{j'}$.
This allows us to specify several cases when the
differential equation can be simplified. 
As an example we consider three-dimensional unimodular Lie groups.
Note that for four-dimensional unimodular Lie groups the Ricci flow
equation was considered in details in \cite{Isenberg}.

\section*{Part I: Ricci Flow of left invariant metrics} 


Let $G$ be a Lie group. To deal with the Ricci flow equation for left
invariant metrics, 
it is convenient to consider the Lie
algebra  $\g$ of $G$ as the Lie algebra of left invariant vector
fields on $G$. Then any left invariant metric tensor field $g$ on $G$
is uniquely determined by the scalar product  
$$
\langle\ ,\  \rangle : \g \times \g \to \R
$$
as follows:
$$
\forall X,Y \in \g, 
\quad g_p(X(p),Y(p)) = 
\langle X, Y\rangle \text{ for any  } p \in G.
$$

Let $\nabla$ be the Riemannian connection of $g$. Then for any
left invariant frame field $E_i$, the connection coefficients
$\Gamma^k_{ij}$ determined by  
$$
\nabla_{E_i} E_j = \Gamma^k_{ij} E_k
$$
are constant.

Now, for the reader's convenience, let us recall calculation of  the
coordinates of Ricci tensor of $g$ with respect to
the left invariant orthonormal frame field.

\subsection*{Ricci tensor of left invariant metric}

\begin{prop}
Let $\{E_i\}$ be an orthonormal frame  of $\g$ with respect to 
$\langle \ , \ \rangle$. Then  
\begin{equation}\label{conex-alg}
\Gamma_{ij}^k=\fracl{1}{2}[C_{ij}^k + C_{ki}^j + C_{kj}^i ],
\end{equation}
where  $[E_i,E_j] = c^k_{ij}E_k$.
\end{prop}

\begin{proof}
Differentiating $\langle E_i, E_j \rangle =\delta_{ij}$, we have
\[
\langle \nabla_{E_k}E_i, E_j \rangle + \langle E_i , \nabla_{E_k}
E_j \rangle =0 \Rightarrow \Gamma_{ki}^j + \Gamma_{kj}^i = 0.
\]
On the other hand the structure constants and the torsion tensor are
defined by
\[
[E_i, E_j]=C_{ij}^k E_k, \qquad T_{ij}^k E_k =
T(E_i,E_j)=\nabla_{E_i}E_j-\nabla_{E_j}E_i-[E_i,E_j].
\]
If the metric is  torsion-free, then we have
\[
\Gamma_{ij}^k - \Gamma_{ji}^k = C_{ij}^k.
\]
From the last equation, by permuting $\{i,j,k\}$, we obtain 
\[
(a)\qquad \Gamma_{ij}^k - \Gamma_{ji}^k = C_{ij}^k,
\]
\[
(b)\qquad \Gamma_{ki}^j - \Gamma_{ik}^j = C_{ki}^j,
\]
\[
(c)\qquad \Gamma_{jk}^i - \Gamma_{kj}^i = C_{jk}^i.
\]
Then $(a)+(b)-(c)$ gives 
\begin{equation}\label{gamma-eq}
\Gamma_{ij}^k=\fracl{1}{2}[C_{ij}^k + C_{ki}^j - C_{jk}^i]
=\fracl{1}{2}[C_{ij}^k  + C_{ki}^j + C_{kj}^i ].
\end{equation}
\end{proof}

\begin{prop}
Let $\{E_i\}$ be an orthonormal frame of $\g$, and $[E_i,E_j] = c^k_{ij}E_k$.
The coordinates  of Ricci tensor with respect to an orthonormal frame
$\{E_i\}$ of $\g$ with respect to $\langle \ , \  \rangle$ 
are 
\begin{equation}
  R_{jk} = \overset{1}{R}_{jk} + \overset{2}{R}_{jk} +
\overset{3}{R}_{jk} + \overset{4}{R}_{jk},
\label{eq:1}
\end{equation}
where
\begin{equation}
  \label{eq:2}
  \overset{1}{R}_{jk} = -\fracl{1}{2} c^s_{mj} c^m_{sk},
\end{equation}
\begin{equation}
  \label{eq:3}
  \overset{2}{R}_{jk}=\fracl{1}{2} \sum_{s=1,n}^n c^m_{ms}(c^k_{sj} + c^j_{sk}),
\end{equation}
\begin{equation}
  \label{eq:4}
  \overset{3}{R}_{jk} = \fracl{1}{4}\sum_{s,m=1}^n c^j_{sm} c^k_{sm},
\end{equation}
\begin{equation}
  \label{eq:5}
  \overset{4}{R}_{jk} = -\fracl{1}{2}\sum_{s,m=1}^n c^m_{sj} c^m_{sk}.
\end{equation}
\end{prop}

\textbf{Proof.}  From the definition of the Riemann curvature tensor
we have
\[
R_{ijk}^{l}E_l \doteq \nabla_{E_i}\nabla_{E_j}E_k -
\nabla_{E_j}\nabla_{E_i} E_k -\nabla_{[E_i, E_j]} E_k,
\]
hence 
\[
R_{ijk}^{l}E_l =\Gamma_{jk}^{l}\Gamma_{il}^{s} E_s -
\Gamma_{ik}^{l}\Gamma_{jl}^{s} E_s
 - C_{ij}^{l}\Gamma_{lk}^{s}E_s.
\]
This implies that
\[
R_{ijk}^{s} =\Gamma_{jk}^{l}\Gamma_{il}^{s} -
\Gamma_{ik}^{l}\Gamma_{jl}^{s}
 - C_{ij}^{l}\Gamma_{lk}^{s},
\]
then we obtain the Ricci tensor:
\begin{equation}
R_{jk} = R_{sjk}^{s} =\Gamma_{jk}^{l}\Gamma_{sl}^{s} -
\Gamma_{sk}^{l}\Gamma_{jl}^{s}
 - C_{sj}^{l}\Gamma_{lk}^{s}.
\label{ricci-formula}
 \end{equation}
From this, using   (\ref{gamma-eq}),   we get
\begin{equation}
R_{jk}= -\fracl{1}{2}C^s_{mj}C^m_{sk} +\fracl{1}{4}C^j_{ms}C^k_{ms}
-\fracl{1}{2}C^s_{mj}C^s_{mk} + \fracl{1}{2}C^m_{ms}(C^k_{sj}+C^j_{sk}).
\end{equation}

\begin{rem}
  $\overset{1}{R}_{ij}$ are coordinates of a tensor (in fact
  $\overset{1}{R}_{ij}$ is the
  Ricci tensor of the Cartan connection on $G$), and the other
  $\overset{\alpha}{R}_{ij}$,
$\alpha=2,3,4$ are not.
\end{rem}

\subsection*{Ricci flow equation for left invariant metrics}

We will rewrite the Ricci flow equation (\ref{rf}) in terms of the 
orthonormal frames. For a smooth one-parameter family $g(t)$  of left invariant metrics  on
$G$, we take a smooth one-parameter family  $\{E_i(t)\}$ such that
$g(t)(E_i(t),E_j(t)) = \delta_{ij}$ (one can easily prove that for any
orthonormal frame $\{E_i\}$ with respect to $g_0$ we can find $\{E_i(t)\}$
such that $E_i(0) = E_i$). 

Now, given a left invariant metric $g_0$,  fix an orthonormal frame
$\{E_i\}$ of $g_0$, then $\{E_i(t)\}$ is uniquely  determined by a 
curve $Q(t)$ in $GL(n)$, where
\begin{equation}
E_i(t) = Q^j_i(t) E_j.
\label{eq:2-1}
\end{equation}

Then the structure constants of $\g$ with respect to $E_i(t)$ are expressed
by 
\begin{equation}
C_{ij}^k(t)=Q_i^s(t)Q_j^m(t)C_{sm}^l\wt Q_l^k(t),
\end{equation}
where $||\wt Q^i_j|| = ||Q^i_j||^{-1}$, 
and for the coordinates of the Ricci tensor $R(t)$ of $g(t)$ we have
\begin{equation}\label{ricci-eq}
  \begin{split}
    R_{jk}(t)= 
 -\fracl{1}{2}C^s_{mj}(t)C^m_{sk}(t)
    +\fracl{1}{4}C^j_{ms}(t)C^k_{ms}(t)
    -\fracl{1}{2}C^s_{mj}(t)C^s_{mk}(t)
    +
\\
\fracl{1}{2}C^m_{ms}(t)[C^k_{sj}(t)+C^j_{sk}(t)].
\end{split}
\end{equation}

\begin{prop} \label{prop:3}
The family $g(t)$ is a solution of the Ricci flow equation \emph{ (\ref{rf})} if and
only if 
 the curve $Q(t)$ defined by \emph{ (\ref{eq:2-1})} satisfies
\begin{equation}
  \label{rf-homo}
  \wt Q^j_s \dot Q^s_i +   \wt Q^i_s \dot Q^s_j = 2 R_{ij}(t).
\end{equation}

\end{prop}
\begin{proof}
From  $g(t)(E_i(t),E_j(t))=\delta_{ij}$ we get 
$$
(\fracl{d}{dt} g(t))(E_i(t),E_j(t)) + g(t)(\fracl{d}{dt} E_i(t),E_j(t))
+ g(t)(E_i(t),\fracl{d}{dt}E_j(t)) = 0.
$$
Then, we have $\fracl{d}{dt} g(t) = -2R$, 
and 
$$
\fracl{d}{dt} E_i(t) = (\fracl{d}{dt} Q_i^s(t)) E_s(0) = 
(\fracl{d}{dt} Q_i^s(t)) \wt Q_s^j E_j(t),
$$
hence follows  (\ref{rf-homo}).
\end{proof}

\subsection*{Reduction to the subgroup of upper triangular  matrices}

From proposition \ref{prop:3} it follows that in order  
to solve the Ricci flow equation (\ref{rf})
we need to find a curve $Q(t)$ in $GL(n)$ satisfying (\ref{rf-homo}).
The following proposition demonstrates that it is sufficient to take
$Q(t)$ in the subgroup $UT(n)$ of upper triangular matrices.

Let $Q(t)$ be a solution of (\ref{rf-homo}). Using the
QR-decomposition,  one can demonstrate that 
\begin{equation}
  Q(t)=B(t)U(t),
\label{eq:QR}
\end{equation}
where 
$B(t)$ is a smooth curve in $UT(n)$ and $U(t)$ is a smooth curve in
the group $O(n)$ of orthogonal matrices.

\begin{prop} \label{prop:4}
A curve   $Q(t)$ in $GL(n)$ is a solution to \emph{(\ref{rf-homo}}) if and
only if the curve $B(t)$ in $UT(n)$ determined by (\emph{\ref{eq:QR}})
is a solution to \emph{(\ref{rf-homo})}.
\end{prop}
\begin{proof}
We take the vector spaces
$\Lambda^2_1 = \{C^k_{ij} \mid C^k_{ij} = -C^k_{ji}\}$ and 
$S^2 = \{A_{ij} \mid A_{ij} = A_{ji}\}$,
and 
  consider  the standard right $GL(n)$-representations:
$$
(T_\Lambda(Q)(C_{uv}^w))^k_{ij} = Q_i^s(t)Q_j^m(t)C_{sm}^l\wt Q_l^k, 
\text{ and } 
(T_S(Q)(A_{uv}))_{ij} = Q_i^p Q_j^q A_{pq}.
$$

We have maps $\overset{\alpha}{R} : \Lambda^2_1 \to S^2$, $\alpha=\overline{1,4}$,
given by (\ref{eq:2})--(\ref{eq:5}).

\begin{lm} \label{lemma:1}
\emph{1)} For any $Q \in GL(n)$, we have
 $\overset{1}{R}(T_\Lambda(Q)C) =  T_S(Q)\overset{1}{R}(C)$ and 
 $R(T_\Lambda(Q)C) =  T_S(Q)R(C)$.

\emph{2)} For any $U \in O(n)$, we have
 $\overset{\alpha}{R}(T_\Lambda(U)C) =  T_S(U)\overset{\alpha}{R}(C)$,
 and thus $R(T_\Lambda(U)C) =  T_S(U)R(C)$.
\end{lm}
\begin{proof}
Direct calculation.  
\end{proof}

Now we can write (\ref{rf-homo}) as 
$$
Q^{-1}(t)\fracl{d}{dt}Q(t) + {}^T(Q^{-1}(t)\fracl{d}{dt}Q(t)) =
2R(T_\Lambda(Q(t))C).
$$
If $Q(t) = B(t)U(t)$, then 
$$
Q^{-1}(t)\fracl{d}{dt}Q(t) = U^{-1}(t)B(t)^{-1}(\fracl{d}{dt}B(t))U(t) + 
U^{-1}(t)\fracl{d}{dt}U(t).
$$
By definition,  $U(t)$ is a curve in $O(n)$, therefore $U^{-1}(t) = {}^TU(t)$,
and $U^{-1}\fracl{d}{dt}U(t)$ lies in the Lie algebra $\mathfrak{o}(n)$ of the Lie
group $O(n)$, hence is a skewsymmetric matrix. Thus we have 
\begin{multline*}
  Q^{-1}(t)\fracl{d}{dt}Q(t) + {}^T(Q^{-1}(t)\fracl{d}{dt}Q(t)) = \\
  {}^TU(t)[B^{-1}(t)\fracl{d}{dt}B(t) +
  {}^T(B^{-1}(t)\fracl{d}{dt}B(t))]U(t)
\end{multline*}
By Lemma (\ref{lemma:1}), 
\begin{multline*}
  R(T_\Lambda(Q(t))C) = R(T_\Lambda(B(t)U(t))C) =
\\
  R(T_\Lambda(U(t))T_\Lambda(B(t))C) =
  {}^TU(t)R(T_\Lambda(B(t))C)U(t),
\end{multline*}
Thus, $Q(t)=B(t)U(t)$  is a solution to (\ref{rf-homo}) if and only if
$Q(t)=B(t)U(t)$  is a solution to (\ref{rf-homo}).
\end{proof}

\begin{rem}
  From proposition \ref{prop:4} it follows that we can take $Q(t)$ in
  the Lie group $UT(n)$, which is smaller than $GL(n)$. However, the
  equation system (\ref{rf-homo}) remains huge, so it is difficult
  even to write down it explicitely, to say nothing about solution.
At the same time, we know that $\overset{1}{R}$ is a symmetric 2-form,
so it is quite natural to choose the initial orthonormal frame $\{E_i(0)\}$
such that $\overset{1}{R}$ has a diagonal matrix with respect to it.
We may suppose that   (\ref{rf-homo}) simplifies in this case.
In the next section we exemplify  this idea in case $G$ is a
three-dimensional unimodular group. 
\end{rem}

\section*{Part II: Three-dimensional unimodular groups}

Let $G$ be a three-dimensional unimodular Lie group. Since the
unimodular group has the property that $C^s_{sk}=0$, we have 
$\overset{2}{R} = 0$. 

\begin{lm}
A three-dimensional Lie group $G$ is unimodular if and only if 
with respect to any frame $E_i$ the Lie algebra $\g$ of $G$ has
structure constants: 
\begin{equation*}
\begin{array}{c} 
C_{ij}^1=
\left(%
\begin{array}{ccc}
  0 & b_1 & -b_2 \\
  -b_1 & 0 & a_1 \\
  b_2 & -a_1 & 0 \\
\end{array}%
\right)\!, \quad C_{ij}^2=
\left(%
\begin{array}{ccc}
  0 & b_3 & -a_2 \\
  -b_3 & 0 & b_2 \\
  a_2 & -b_2 & 0 \\
\end{array}%
\right)\!, 
\\[25pt] 
C_{ij}^3=
\left(%
\begin{array}{ccc}
  0 & a_3 & -b_3 \\
  -a_3 & 0 & b_1 \\
  b_3 & -b_1 & 0 \\
\end{array}%
\right)\!,
\end{array}
\end{equation*}

and 

\[
\overset{1}R_{jk}=\left(%
\begin{array}{ccc}
  -2a_2a_3+2b_3^2 & 2a_3b_2-2b_1b_3 & 2a_2b_1-2b_2b_3 \\
  2a_3b_2-2b_1b_3 & -2a_1a_3+2b_1^2 & -2b_1b_2+2a_1b_3 \\
  2a_2b_1-2b_2b_3 & -2b_1b_2+2a_1b_3 & -2a_1a_2+2b_2^2 \\
\end{array}%
\right).
\]

In particular, the matrix $||\overset{1}R_{jk}||$ is diagonal if and
only if 
\begin{equation}
  \label{eq:24}
    b_1b_3-a_3b_2=0, \quad 
    b_2b_3-a_2b_1=0, \quad
    b_1b_2-a_1b_3=0.
\end{equation}
\end{lm}

\begin{proof}
The set of structure constants $C^k_{ij}$ given above is the general
solution of the equation system 
$C^k_{ij} = -C^k_{ji}$, $C^s_{sm} = 0$. Note that in the
three-dimensional case  these equations imply the Jacobi identities.
\end{proof}

Solving  the equation system (\ref{eq:24}), we obtain the following
cases:

\subsection*{Case I: $b_1=b_2=b_3=0$}

The  structure constants are
\[
C_{ij}^1=
\left(%
\begin{array}{ccc}
  0 & 0 & 0 \\
  0 & 0 & a_1 \\
  0 & -a_1 & 0 \\
\end{array}%
\right)\!,
 \quad
C_{ij}^2=
\left(%
\begin{array}{ccc}
  0 & 0 & -a_2 \\
  0 & 0 & 0 \\
  a_2 & 0 & 0 \\
\end{array}%
\right)\!,
 \quad
 C_{ij}^3=
\left(%
\begin{array}{ccc}
  0 & a_3 & 0 \\
  -a_3 & 0 & 0 \\
  0 & 0 & 0 \\
\end{array}%
\right)\!.
\]

\begin{rem}
Note that, if  $a_1=a_2=a_3=1$, the algebra $\mathfrak{g}\cong \frak{so}(3)$.
\end{rem}
The nonzero coordinates of  Ricci tensor are:
\[
\begin{array}{ccc}
R_{11}=  \fracl{1}{2}(a_1^2-(a_2-a_3)^2) \\
R_{22}= \fracl{1}{2}(-a_1^2+a_2^2+2a_1a_3-a_3^2) \\
R_{33}=  \fracl{1}{2}(-a_1^2+2a_1a_2-a_2^2+a_3^2) \\
\end{array}%
\]
and 
\[
\overset{1}R_{jk}(t)=\left(%
\begin{array}{ccc}
  a_2a_3 & 0 & 0 \\
  0 & a_1a_3 & 0 \\
   0 & 0 & a_1a_2 \\
 \end{array}%
 \right)\!.
 \]

Thus with respect to this frame both tensors $R$ and $\overset{1}{R}$
have diagonal form.

If we take  $Q(t) \in \Delta(3)$, where $\Delta(n)$ is the
subgroup of nondegenerate diagonal matrices in $GL(n)$, such that
\[
Q(t)=\left(%
\begin{array}{ccc}
  f(t) & 0 & 0 \\
  0 & g(t) & 0 \\
  0 & 0 & h(t) \\
\end{array}%
\right)\!,
\]
then the nonzero coordinates of Ricci tensor are:
\begin{equation} \label{ricci1}
\begin{array}{l}
R_{11}(t)=
a_2a_3f^2(t)+\fracl{a_1^2g^2(t)h^2(t)}{2f^2(t)}-\fracl{f^2(t)(a_3^2g^4(t)+a_2^2h^4(t))}{2g^2(t)h^2(t)},
\\
R_{22}(t)=
a_1a_3g^2(t)+\fracl{a_2^2f^2(t)h^2(t)}{2g^2(t)}-\fracl{g^2(t)(a_3^2f^4(t)+a_1^2h^4(t))}{2f^2(t)h^2(t)},
\\
R_{33}(t)=\fracl{a_3^2f^4(t)g^4(t)-(a_2f^2(t)-a_1g^2(t))^2h^4(t)}{2f^2(t)g^2(t)h^2(t)}.
\end{array}
\end{equation}

Then the Ricci flow equation system is 
\begin{equation}
\label{case-i}
  \begin{array}{l}
\fracl{1}{f(t)} \fracl{df}{dt}(t) = 
R_{11}(t),
\\
\fracl{1}{g(t)} \fracl{dg}{dt}(t) = 
R_{22}(t),
\\
\fracl{1}{h(t)} \fracl{dh}{dt}(t) =
R_{33}(t),
  \end{array}
\end{equation}
where $R_{11}(t)$,  $R_{22}(t)$, and  $R_{33}(t)$ are given by (\ref{ricci1}).
\emph{Thus, in \emph{Case I}  the Ricci flow ODE system
  \emph{(\ref{rf-homo})} reduces to the ODE system
\emph{(\ref{case-i})} 
in three unknown functions.}

\subsection*{Case II:  $b_1  \ne 0$,  $b_2=b_3=0$}

The  structure constants are
\[
\begin{array}{c}
C_{ij}^1=
\left(%
\begin{array}{ccc}
  0 & b_1 & 0 \\
  -b_1 & 0 & a_1 \\
  0 & -a_1 & 0 \\
\end{array}%
\right)\!,
 \quad 
C_{ij}^2=
\left(%
\begin{array}{ccc}
  0 & 0 & 0 \\
  0 & 0 & 0 \\
  0 & 0 & 0 \\
\end{array}%
\right)\!,
\\[25pt]
  C_{ij}^3=
\left(%
\begin{array}{ccc}
  0 & a_3 & 0 \\
  -a_3 & 0 & b_1 \\
  0 & -b_1 & 0 \\
\end{array}%
\right)\!.
\end{array}
\]

 The Ricci tensor has coordinates
\[
R_{jk}(t)=\left(%
\begin{array}{ccc}
  \fracl{1}{2}(a_1^2-a_3^2) & 0 & (a_1+a_3)b_1 \\
  0 & \fracl{1}{2}(-a_1^2+2a_1a_3-a_3^2-4b_1^2) & 0 \\
  (a_1+a_3)b_1 & 0 & \fracl{1}{2}(-a_1^2+a_3^2) \\
\end{array}%
\right),
\]
and 
\[
\overset{1}R_{jk}(t)=\left(%
\begin{array}{ccc}
  0 & 0 & 0 \\
  0 & -a_1a_3+b_1^2 & 0 \\
  0 & 0 & 0 \\
\end{array}%
\right)\!.
\]

If we take  $Q(t) \in UT(3)$ such that
\[
Q(t)=\left(%
\begin{array}{ccc}
  f(t) & 0 & 0 \\
  0 & g(t) & 0 \\
  w(t) & 0 & h(t) \\
\end{array}%
\right),
\]
then the nonzero coordinates of Ricci tensor are:
\begin{equation} \label{ricci2}
  \begin{array}{l}
R_{11}(t) = 
- g^2(t) [a_1^2  w^4(t) - 4 a_1 b_1 f(t)  w^3(t) 
+
\\[2pt] 
\quad (4 b_1^2 + 2 a_1 a_3) f^2(t) w^2(t) - 
4 a_3 b_1 f^3(t)  w(t)
- 
\\[2pt]
\quad a_1^2 h^4(t) + a_3^2 f^4(t)]/(2 f^2(t) h^2(t)) 
\\[5pt]
R_{13}(t) =
-g^2(t) [a_1^2  w^3(t) - 3 a_1 b_1 f(t)  w^2(t) 
+ 
\\[2pt]
\quad (a_1^2 h^2(t) + 
(2 b_1^2 + a_1^2 a_3) f^2(t)) w(t) 
- 
\\[2pt]
\quad a_1 b_1 f(t)h^2(t) - a_3 b_1 f^3(t)]/(f^2(t) h(t))
\\[5pt]
R_{22}(t)= 
- g^2(t) [a_1^2  w^4(t) - 4 a_1 b_1 f(t)  w^3(t)
+ 
\\[2pt]
\quad (2 a_1^2  h^2(t) + 
(4 b_1^2 + 2 a_1 a_3) f^2(t)) w^2(t) -
\\[2pt]
\quad 4 (a_1 b_1 f(t) h^2(t) + a_3 b_1 f^3(t)) w(t) +
a_1^2 h^4(t)
+ 
\\[2pt]
\quad (4 b_1^2 - 2 a_1 a_3) f^2(t) h^2(t)  + a_3^2 f^4(t)]
/(2f^2(t) h^2(t))
\\[5pt]
R_{33}(t)= 
g^2(t)[a_1^2  w^4(t) - 4 a_1 b_1 f(t)  w^3(t)
+ 
\\[2pt]
\quad (4 b_1^2 + 2 a_1 a_3) f^2(t)  w^2(t) 
- 4 a_3 b_1 f^3(t) w(t) 
-
\\[2pt]
\quad a_1^2 h^4(t)  +  a_3^2 f^4(t)]/(2 f^2(t) h^2(t)) 
\end{array}
\end{equation}

Then the Ricci flow equation system is 
\begin{equation} 
\label{case-ii}
  \begin{array}{l}
\fracl{1}{f(t)} \fracl{df}{dt}(t) = R_{11}(t) 
\\[8pt]
\fracl{1}{g(t)} \fracl{dg}{dt}(t) = 
R_{22}(t)
\\[8pt]
\fracl{1}{h(t)} \fracl{dh}{dt}(t) =
R_{33}(t)
\\[8pt]
{}[f(t)\fracl{d}{dt}w(t) - w(t)\fracl{d}{dt}f]/(f(t)h(t)) = 2 R_{13}(t),
  \end{array}
\end{equation}
where $R_{11}(t)$, $R_{22}(t)$, $R_{33}(t)$, and  $R_{13}(t)$ are
given by (\ref{ricci2}).

\emph{Thus, in \emph{Case II}  the Ricci flow ODE system
  \emph{(\ref{rf-homo})} reduces to the ODE system
\emph{(\ref{case-ii})} 
in four unknown functions.}

\subsection*{Case III: $b_1 \ne 0$, $b_2 \ne 0$, $b_3 \ne 0$}

We can find $\rho$,$\alpha$, and $\beta$ such that 
$$
\begin{array}{l}
b_1 = \sqrt{\rho \fracl{\cos(\beta) \cos(\alpha)}{\sin(\beta)}}
\\
b_2 = \sqrt{\rho \fracl{\sin(\beta) \cos(\alpha)}{\cos(\beta)}}
\\
b_3 = \sin(\alpha)\sqrt{\rho \fracl{\sin(\beta) \cos(\beta)}{\cos(\alpha)}}
\end{array}
$$

Now we take the new frame 
$$
\begin{array}{l}
E_1 = (\cos(\alpha), \sin(\alpha)\sin(\beta), \sin(\alpha)\cos(\beta))
\\
E_2 = (\sin(\alpha), -\cos(\alpha)\sin(\beta), -\cos(\alpha)\cos(\beta))
\\
E_3 = (0, \cos(\beta), -\sin(\beta))
\end{array}
$$
With respect to this frame the only nonzero $C^k_{ij}$ are 
$$
C^1_{23} = - C^1_{32} = 
\fracl{1}{\sin(\alpha)}\sqrt{\fracl{\rho}{\cos(\alpha)\cos(\beta)\sin(\beta)}}
$$
and\emph{ we arrive at} Case I.

\begin{rem}
One can easily see that the case $b_1 \ne 0$, $b_2 \ne 0$, and $b_3 =0$ is impossible.  
\end{rem}

\begin{rem}
The calculations were performed by computer system Maxima.  
\end{rem}



\begin{thebibliography}{99}

\bibitem{Besse}
Besse A.L., \newblock{\em Einstein Manifolds}, Springer-Verlag,
Berlin (1987)

\bibitem{Cao-Chow}
Cao H-D., Chow B., \newblock{Recent developements on the Ricci
flow},
\newblock{\em Bull.(N.S.) Amer.Math.Soc.}\textbf{36} (1999), 59-74.

\bibitem{Cartan}
Cartan \'E., \newblock{\em La M\'etode du rep\'ere mobile, La th\'orie des
groupes continus et les espaces g\'en\`eralis\'e}, Conf\'erences faites \`a
Moscou, 16-20 Juin 1930.  \newblock{\em Gosudarsvennoe
Technicheskoe-Teoreticheskoe Isdatelstvo}, Moscow, (1933)

\bibitem{Chow}
Chow B., \newblock{ A survey of Hamilton's Program for the Ricci
flow on 3-manifolds}, arXiv:math.DG/0211266v1 18 Nov. 2002.

\bibitem{Dubrovin-Novikov-Fomenko}
Dubrovin B.A., Novikov S.P., Fomenko A.T. \newblock{\em Modern
Geometry}, Nauka, Moscow  (1979)

\bibitem{Hamilton}
Hamilton R.S.,
\newblock {Three-manifolds with positive Ricci curvature}.
\newblock {\em J.Differential Geom. }\textbf{17} (1982), no. 2, 255-306

\bibitem{Isenberg}
Isenberg J., Jackson M., Peng L., \newblock{ Ricci flow on locally
homogeneous closed 4-manifolds}, 
arXiv:math.DG/0502170v1 8 Feb. 2005.


\bibitem{Kobayashi-Nomizu}
Kobayashi S., Nomizu K., \newblock{\em Foundations of Differential
geometry}, John Wiley \& Sons, Inc., Vol II, New York (1979)

\bibitem{Milnor}
Milnor J., \newblock{Towards the Poincar\'e Conjecture and the
Classification of $3$-Manifolds}, \newblock{\em Notices of the AMS},
Vol. \textbf{50}, Number \textbf{10}, 1226-1233 (2003)

\bibitem{Morgan}
Morgan J.W., \newblock{ Recent Progress on the Poincar\'e Conjecture
and the Classification of 3-manifolds} \newblock{\em Bull (N.S.)
Amer.Math.Soc.} Vol \textbf{42} Num. \textbf{1} (2004) 57-78

\bibitem{Perelman}
Perelman G., \newblock{ Ricci flow with surgery on three-manifolds},
arXiv:math.DG/0303109v1 10 Mar 2003.

\bibitem{Sesun}
Sesun N., \newblock{ Convergence of the Ricci flow toward a unique
soliton},  arXiv:math.DG/0405398v1 20
May 2004.

\bibitem{Shapukov}
Shapukov B.N., \newblock{\em Grupos y Algebras de Lie en ejercicios
y problemas}, Matematika URSS (2001)

\end{thebibliography}
\end{document}